\documentclass[a4paper,11pt]{article}
\usepackage{amscd, amsmath, amssymb, amsfonts, epic}
\usepackage{graphicx}
\usepackage{a4wide}
\usepackage{xypic}
\usepackage{paralist}
\usepackage{booktabs}
\usepackage{multirow}

\title{The cone conjecture for some rational elliptic threefolds}
\author{Arthur Prendergast-Smith}

\date{}

\def\Z{\text{\bf Z}}
\def\Q{\text{\bf Q}}
\def\R{\text{\bf R}}

\def\P{\text{\bf P}}

\def\arrow{\rightarrow}

\def\iso{\cong}
\def\curlyo{\mathcal{O}}
\def\Pic{\text{Pic}} 
\def\Aut{\text{Aut}}
\def\PsAut{\text{PsAut}}
\def\Div{\text{Div}}

\newcommand{\Nef}[1]{\overline{A(#1)}}
\newcommand{\Mov}[1]{\overline{M(#1)}}
\newcommand{\Curv}[1]{\overline{\text{Curv}(#1)}}

\newtheorem{theorem}{Theorem}[section]
\newtheorem{lemma}[theorem]{Lemma}
\newtheorem{corollary}[theorem]{Corollary} 
\newtheorem{proposition}[theorem]{Proposition}

\newtheorem{conjecture}[theorem]{Conjecture}

\begin{document}

\maketitle

A central problem of modern minimal model theory is to describe the
various cones of divisors associated to a projective variety. For Fano
varieties the nef cone and movable cone are rational polyhedral by the
cone theorem \cite[Theorem 3.7]{KollarMori1998} and the theorem of
Birkar--Cascini--Hacon--McKernan \cite{BCHM2007}. For more general
varieties the picture is much less clear: these cones need not be
rational polyhedral, and can even have uncountably many extremal rays.

The Morrison-Kawamata cone conjecture \cite{Morrison1992,
  Kawamata1997, Totaro2008} describes the action of automorphisms on
the cone of nef divisors and the action of pseudo-automorphisms on the
cone of movable divisors, in the case of a Calabi-Yau variety, a
Calabi-Yau fibre space, or a Calabi-Yau pair. Although these cones
need not be rational polyhedral, the conjecture predicts that they
should have a rational polyhedral fundamental domain for the action of
the appropriate group. It is not clear where these automorphisms or
pseudo-automorphisms should come from; nevertheless, the conjecture
has been proved in various contexts by Sterk--Looijenga--Namikawa
\cite{Sterk1985,Namikawa1985} Kawamata \cite{Kawamata1997}, and Totaro
\cite{Totaro2008b}.

In this paper we give some new evidence for the conjecture, by
verifying it for some threefolds which are blowups of $\P^3$ in the
base locus of a net (that is, a 2-dimensional linear system) of
quadrics. Our main result is the following:

\begin{theorem} \label{theorem-mainthm}
  Let $X$ be the blowup of $\P^3$ in $8$ distinct points which are the
  base locus of a net of quadrics.

  \quad (1) The nef cone $\Nef{X}$ is rational polyhedral and spanned by
  effective divisors.

  \quad (2) If the net has no reducible member, the effective movable
  cone $\Mov{X}^e$ has a rational polyhedral fundamental domain for
  the action of $PsAut(X)$.
\end{theorem}

In Section \ref{sect-prelims} we will see that for any net of quadrics
in $\P^3$ with $8$ distinct basepoints, the blowup of the base locus
of the net has an elliptic fibration over $\P^2$. The condition that
the net have no reducible member is equivalent to the generic fibre of
the fibration (an elliptic curve over the function field of $\P^2$)
having Mordell--Weil rank 7, the maximum possible. Although in
statement (2) we restrict to this class of nets, much of the proof
works for general nets, and it should be possible to fill in the
remaining details.

The proof of Theorem \ref{theorem-mainthm} relies to a large extent on
the explicit geometry of nets of quadrics in $\P^3$ and so seems
difficult to generalise to other classes of varieties. Nevertheless,
the result is significant inasmuch as it seems to be the first
verification of the cone conjecture for a klt Calabi--Yau pair
$(X,\Delta)$ of dimension 3 with $\Delta \neq 0$. (See the next
section for definitions.) This lends further support to the point of
view that klt Calabi--Yau pairs provide a natural setting for the
conjecture.

Thanks to Klaus Hulek and Burt Totaro for their comments.

\section{The cone conjecture} \label{section-conjecture}
In this section we give the precise statement of the cone conjecture
for klt Calabi--Yau pairs, following \cite{Totaro2008b}. (See also
Section 1 of \cite{Totaro2008b} for history and examples.) We work
throughout over an algebraically closed field $k$ of characteristic
0. A {\it rational polyhedral cone} in a real vector space $V$ with a
$\Q$-structure is a closed convex cone with finitely many extremal
rays, each spanned by a rational vector.

Suppose $f: X \arrow S$ is a projective surjective morphism of normal
varieties with connected fibres. A Cartier divisor $D$ on $X$ is said
to be {\it $f$-nef} (resp. {\it $f$-movable}, {\it $f$-effective}) if
$D \cdot C \geq 0$ for all curves $C$ mapped to a point by $f$ (resp.
if $\operatorname{codim \ Supp \ Coker } (f^*f_* \mathcal{O}_X(D)
\arrow \mathcal{O}_X(D)) \geq 2$, if $f_* \mathcal{O}_X(D) \neq 0$).

We define the real vector space $N^1(X/S)$ to be $\Div(X)/\cong_S
\otimes \ \R$ where $\Div(X)$ is the group of Cartier divisors on $X$
and $\cong_S$ denotes numerical equivalence over $S$. We denote by
$N^1(X/S)_\Z$ the free abelian group in $N^1(X/S)$ consisting of
numerical classes of Cartier divisors. The {\it $f$-nef} cone
$\Nef{X/S}$ (resp. {\it closed $f$-movable cone} $M(X/S)$, {\it
  $f$-pseudoeffective cone} $B(X/S)$) is the closed convex cone
generated by classes of {\it $f$-nef} (resp. {\it $f$-movable}, {\it
  $f$-effective}) divisors. The {\it $f$-effective cone} $B^e(X/S)$ is
the cone generated by $f$-effective Cartier divisors. We denote by
$\Nef{X/S}^e$ and $\Mov{X/S}^e$ the intersections $\Nef{X/S} \cap
B^e(X/S)$ and $M(X/S) \cap B^e(X/S)$, and call them the {\it
  $f$-effective $f$-nef cone} and {\it $f$-effective $f$-movable cone}
respectively.

Define a {\it pseudo-isomorphism} from $X_1$ to $X_2$ over $S$ to be a
birational map $X_1 \dashrightarrow X_2$ over $S$ which is an
isomorphism in codimension 1. A {\it small $\Q$-factorial
  modification} (SQM) of $X$ over $S$ means a pseudo-isomorphism over $S$
from $X$ to another $Q$-factorial variety with a projective morphism
to $S$.

For an $\R$-divisor $\Delta$ on a normal $\Q$-factorial variety $X$,
the pair $(X,\Delta)$ is {\it klt} if, for all resolutions $\pi:
\tilde{X} \arrow X$ with a simple normal crossing $\R$-divisor
$\tilde{\Delta}$ such that $K_{\tilde{X}}+\tilde{\Delta} =
\pi^*(K_X+\Delta)$, the coefficients of $\tilde{\Delta}$ are less than
1. (In particular if $X$ is smooth and $D$ is a smooth divisor on $X$,
then $(X,rD)$ is klt for any $r<1$.) We say that $(X/S,\Delta)$ is a
{\it klt Calabi--Yau pair} if $(X,\Delta)$ is a $\Q$-factorial klt
pair with $\Delta$ effective such that $K_X+\Delta$ is numerically
trivial over $S$. 

We denote the groups of automorphisms or pseudo-automorphisms of $X$
over $S$ which preserve a divisor $\Delta$ by $\Aut(X/S,\Delta)$ and
$\PsAut(X/S,\Delta)$. Note that the action of $\Aut(X/S,\Delta)$ and
$\PsAut(X/S,\Delta)$ on $N^1(X/S)$ is determined by the images of the
representations $\Aut(X/S,\Delta) \arrow GL(N^1(X/S)_\Z)$ and
$\PsAut(X/S,\Delta) \arrow GL(N^1(X/S)_\Z)$. We denote the images of
these representations by $\Aut^*(X/S,\Delta)$ and
$\PsAut^*(X/S,\Delta)$.

\begin{conjecture} \label{conj-totaro}
Let $(X/S, \Delta)$ be a klt Calabi--Yau pair. Then:

\quad (1) The number of $\Aut(X/S,\Delta)$-equivalence classes of faces of
the effective nef cone $\Nef{X/S}^e$ corresponding to birational
contractions or fibre space structures is finite. Moreover, there
exists a finite rational polyhedral cone $\Pi$ which is a fundamental
domain for the action of $\Aut^*(X/S,\Delta)$ on $\Nef{X/S}^e$ in the sense that

\quad (a) $\Nef{X/S}^e = \Aut^*(X/S,\Delta) \cdot \Pi$,

\quad (b) $\text{Int} \ \Pi \cap g \text{Int} \ \Pi = \emptyset$ for
$g \neq 1$ in $\Aut^*(X/S,\Delta)$. 

\quad (2) The number of $\PsAut(X/S,\Delta)$-equivalence classes of chambers
$\Nef{X'/S, \alpha}^e$ in the cone $\Mov{X/S}^e$ corresponding to
marked SQMs $f':X' \arrow S$ of $X \arrow S$ with marking $\alpha: X'
\dashrightarrow X$ is finite. Moreover, there exists a finite rational
polyhedral cone $\Pi'$ which is a fundamental domain for the action of
$\PsAut^*(X/S,\Delta)$ on $\Mov{X/S}^e$.
\end{conjecture}

The conjecture has been proved for Calabi--Yau surfaces by
Looijenga--Sterk and Namikawa \cite{Sterk1985,Namikawa1985}, for klt
Calabi--Yau pairs of dimension 2 by Totaro \cite{Totaro2008b}, and for
Calabi--Yau fibre spaces of dimension 3 over a positive-dimensional
base by Kawamata \cite{Kawamata1997}. For Calabi--Yau 3-folds there
are significant results by Oguiso--Peternell \cite{Oguiso2001},
Szendr\"oi \cite{Szendroi1999}, Uehara \cite{Uehara2004}, and Wilson
\cite{Wilson1992}, but the conjecture remains open.

\section{Nets of quadrics in $\P^3$} \label{sect-prelims}

In this section we give some relevant facts about blowups of $\P^3$ in
the base locus of a net of quadrics and fix some notation. We then
explain what the cone conjecture predicts in this situation.

If $X$ is the blowup of $\P^3$ in any set of 8 points ${p_1,\ldots,p_8}$,
then $N^1(X)$ is 9-dimensional with basis $\{H,E_1,\ldots, E_8\}$,
where $H$ is the pullback to $X$ of the hyperplane class on $\P^3$ and
$E_i$ is the class of the exceptional divisor of the blowup of
$p_i$. The dual vector space $N_1(X)$ has basis
$\{l,l_1,\ldots,l_8\}$, where $l$ is the pullback to $X$ of the class
of a line in $\P^3$ and $l_i$ is the class of a line in $E_i$. The
intersection pairing between these spaces is specified by the
following intersection numbers: $H \cdot l=1$, $H \cdot l_i=0$, $E_i
\cdot l=0$, $E_i \cdot E_j = - \delta_{ij}$, for all $i$ and $j$.

Now suppose the 8 points are distinct and are the base locus of a net
of quadrics in $\P^3$. The proper transforms of quadrics in the net
are (up to scalar) sections of the line bundle
$2H-E_1-\ldots-E_8=-\frac12K_X$: since we have blown up the base locus
of the net, $-\frac12K_X$ is basepoint-free on $X$ and so gives a
surjective morphism $f:X \arrow \P^2$. Since $f$ is given by sections
of $-\frac12K_X$ we have $-\frac12K_X=f^*(L)$ for $L$ the hyperplane
class on $\P^2$. This implies that $-\frac12K_X \cdot C=0$ for any
curve $C$ on $X$ mapped to a point by $f$. Adjunction therefore tells
us that the smooth fibres of $f$ are curves with trivial canonical
bundle, hence elliptic curves. In other words, $f: X \arrow \P^2$ is
an elliptic fibration. If $X_\eta$ denotes the generic fibre of $f$,
we define the {\it Mordell--Weil rank} $\rho$ of $f$ (or of $X$) to be
the rank of the finitely-generated abelian group $\Pic^0(X_\eta)$ of
degree-0 line bundles on $X_\eta$. One can show \cite[Theorem
  7.2]{Totaro2008} that $\rho(f)=7-d$ where $d$ is the number of
reducible quadrics (unions of 2 distinct planes) in the net.

The elliptic fibration $f$ on $X$ is important because it gives us a
supply of pseudo-automorphisms of $X$. Using the group law on an
elliptic curve, $\Pic^0(X_\eta)$ acts on $X_\eta$ by automorphisms and
by \cite[Lemma 6.2]{Totaro2008} this extends to an action on $X$ by
pseudo-automorphisms. That is, we can identify $\Pic^0(X_\eta)$ with a
subgroup of $\PsAut(X)$. (We will see in the course of the proof that
this subgroup gives enough pseudo-automorphisms to verify the
conjecture.) More precisely, since $f$ is given by sections of the line
bundle $-\frac12K_X$, the action of elements of $\Pic^0(X_\eta)$
preserves divisors $\Delta=\frac12 D$ for $D$ a smooth divisor in the
linear system $|-2K_X|$, and commutes with the morphism $f:
X\arrow \P^2$, so we can identify $\Pic^0(X_\eta)$ with a subgroup of
$\PsAut(X/\P^2,\Delta)$ for any $\Delta$ of this form.

We must also say something about the reducible fibres of $f$. Note
that by our description of $f$, all fibres are isomorphic to quartic
curves in $\P^3$ which are the complete intersection of 2 quadrics in
the net. First suppose a reducible fibre contains a line $L$. It is
easy to see that $L$ must be the line joining 2 basepoints $p_i$ and
$p_j$ of the net, so its proper transform on $X$ has class
$l-l_i-l_j$ in $N_1(X)$. We denote this class by $C_{ij}$.

We will see in due course that the classes $C_{ij}$ play an important
role in the proof of Theorem \ref{theorem-mainthm}. Note that there
are ${8 \choose 2}=28$ such lines, each contained in exactly 1 fibre
of $f$, and hence at most 28 fibres containing a line.

If a reducible fibre does not contain a line, it is the union of 2
irreducible conics in $\P^3$. Each conic is contained in a plane, and
the union of the planes is a reducible quadric in the net. We will
denote the classes in $N^1(X)$ of the 2 components of the proper
transform of a reducible quadric $Q_i$ in the net by $D^a_i$ ($a=1,\,
2$). For any $i$ we have $D^1_i+D^2_i=-\frac12K_X=f^*(L)$, so both
components must be mapped by $f$ to a line $L_i$ in $\P^2$. For any
point $p \in L_i$ the fibre $f^{-1}(p)$ is then a reducible curve, the
union of 2 conics in $\P^3$, one contained in each component $D^a_i$
of $Q_i$. We denote the class in $N_1(X)$ of the (possibly reducible)
curve $f^{-1}(p) \cap D^a_i$ by $F^a_i$. 

It is easy to see that any such plane and any such conic must both
contain exactly 4 basepoints $p_q,\, p_r,\, p_s, \, p_t$ of the net,
so in terms of our bases for $N^1(X)$ and $N_1(X)$ their proper
transforms have classes $D^a_i=H-E_q-E_r-E_s-E_t$ and
$F^a_i=2l-l_q-l_r-l_s-l_t$. By the intersection numbers given above we
get $D^a_i \cdot F^a_i =-2$. Also if $F$ is the class of any fibre of
$f$ we have $D^a_i \cdot F =0$ because $D^a_i$ maps to a line in
$\P^2$. Since $F=F^1_i+F^2_i$ for any $i$ we get $D^a_i \cdot F^b_i=2$
for $a \neq b$.

Define a prime divisor $D$ on $X$ to be {\it vertical} if $f(D) \neq
\P^2$. Since any divisor pulled back from $\P^2$ is a multiple of
$-\frac12K_X$, the description of the reducible fibres of $f$ shows
that the only vertical divisors on $X$ have divisor class either a
multiple of $-\frac12K_X$ or else $D^a_i$, where the latter are
effective. We will see that vertical divisors play an important role
in describing the movable cone of $X$: namely, Lemma
\ref{lemma-relativemovable} shows that the $f$-movable cone is more or
less defined by intersection numbers with fibral curves lying inside
vertical divisors. Note however that for the final steps of the proof,
we restrict to the case of Mordell--Weil rank 7, which by the
discussion above is equivalent to the fact that $X$ has no vertical
divisors other than multiples of $-\frac12K_X$.

We mention some facts about the birational geometry of
$X$. Suppose $\phi: X \dashrightarrow X'$ is some other projective
variety obtained by flopping some {\it $f$-fibral curves} on $X$ (that
is, curves contained in fibres of $f$). The line bundle
$-\frac12K_{X'}$ is basepoint-free on $X'$ and gives another elliptic
fibration $f': X' \arrow \P^2$ such that $f = f' \circ \phi$ as
rational maps. Also, $\phi$ induces an identification $\phi_*$ of the
spaces $N^1(X)$ and $N^1(X')$ and hence an identification of the dual
spaces $N_1(X)$ and $N_1(X')$. Therefore for any such $X'$ we can
think of the nef cone $\Nef{X'}$ as a cone in the vector space
$N_1(X)$, and the closed cone of curves $\Curv{X'}$ (the dual
of the nef cone) as a cone in $N_1(X)$. Also note that $\phi_*$
identifies $K_X$ and $K_{X'}$ and so the subspaces $K_X^\perp = \{x
\in N_1(X) : K_X \cdot x =0\}$ and $K_{X'}^\perp = \{x \in N_1(X') :
K_X \cdot x =0\}$ are identified. We can therefore speak of the
subspace $K^\perp \subset N_1(X)$ without reference to a particular
model of $X$.

Now we explain the predictions of the cone conjecture in this
situation. If $X$ is the blowup of the base locus of a net of
quadrics, we saw that the line bundle $-\frac12K_X$ is
basepoint-free. Therefore $-2K_X$ is basepoint-free also, so by
Bertini's theorem a general divisor $D \in |-2K_X|$ is smooth. As
mentioned in the previous section the pair $(X,\frac12 D)$ is then
klt, and $K_X+\frac12 D$ is numerically trivial (over $S=\text{ Spec }
k$). The cone conjecture therefore predicts that the groups
$\Aut^*(X,\frac12 D)$ and $\PsAut^*(X, \frac12 D)$ act on the cones
$\Nef{X}^e$ and $\Mov{X}^e$ respectively with rational polyhedral
fundamental domain. The first statement of Theorem
\ref{theorem-mainthm} says that the prediction about the nef cone is
true for all such $X$, in a strong sense: the nef cone itself is
rational polyhedral. (The existence of a rational polyhedral
fundamental domain then follows, as we will see in the next section.)
The second statement of the theorem says that the prediction about the
movable cone is also true, although (as we shall see) that cone itself
is `almost never' rational polyhedral.

\section{Nef cones} \label{sect_nefcones}

In this section we will prove the first statement of Theorem
\ref{theorem-mainthm}, namely that if $X$ is the blowup of $\P^3$ in
the base locus of a net of quadrics with 8 distinct basepoints, then
$\Nef{X}$ is a rational polyhedral cone. In the case where $X$ has
Mordell--Weil rank $\rho=7$, we prove the same thing about the nef
cones of flops $X'$ of $X$ which we will use in the next section. 

 The cone theorem \cite[Theorem
  3.7]{KollarMori1998} says (in any dimension) that if $(X,\Delta)$ is
a klt pair with $\Delta$ effective, any $(K_X+\Delta)$-negative
extremal ray of $\Curv{X}$ can be contracted to give a
projective variety $Z$. In the case that $X$ is s smooth threefold and
$\Delta=0$, the following theorem of Mori \cite[Theorem 3.3, Theorem
  3.5]{Mori1982} gives the possibilities for the exceptional locus of
the contraction:

\begin{theorem}[Mori] \label{theorem-mori3folds}

  Suppose that $X$ is a smooth projective threefold, and $f: X \arrow
  Z$ is the contraction morphism associated to a $K_X$-negative
  extremal ray of $\Curv{X}$. Then either $\operatorname{dim }
  Z \leq 2$ and the anticanonical bundle $-K_X$ is $f$-ample, or else
  $f$ is birational, the exceptional set $\text{Exc}(f)$ is a prime
  divisor $D$ on $X$, and the possibilities for $D$ and $f$ are as
  follows:

\begin{enumerate}
\item $D$ is a $\P^1$-bundle over a smooth curve $C$, and $f_{|D}$ is the
  bundle map $D \arrow C$,
\item $D \iso \P^2$ with normal bundle $\mathcal{O}_D (D) \iso
  \mathcal{O}_{\P^2}(-1)$, and $f$ contracts $D$ to a smooth point,
\item $D \iso \P^1 \times \P^1$ with $\mathcal{O}_D (D)$ of bidegree
  $(-1,-1)$, and $f$ contracts $D$ to a point,
\item $D$ is isomorphic to a singular quadric in $\P^3$ with
  $\mathcal{O}_D (D) = \mathcal{O}_D \otimes \mathcal{O}_{\P^3}(-1)$,
  and $f$ contracts $D$ to a point,
\item $D \iso \P^2$ with normal bundle $\mathcal{O}_D (D) \iso
  \mathcal{O}_{\P^2}(-2)$, and $f$ contracts $D$ to a point.
\end{enumerate}

\end{theorem}

Only one of these possibilities is relevant to us:

\begin{proposition}
  Suppose $X$ is a threefold obtained by blowing up the base locus of
  a net of quadrics in $\P^3$, and let $R$ be a $K_X$-negative
  extremal ray of the closed cone of curves $\Curv{X}$. Then the
  contraction morphism $\operatorname{cont}_R:X \arrow Z$ is
  birational of type $2$ on the above list --- that is, the
  exceptional divisor $D$ is isomorphic to $\P^2$, with normal bundle
  $\mathcal{O}_{\P^2}(-1)$, and $D$ blows down to a smooth
  point. Moreover, the exceptional divisor $D$ in this case must be
  the exceptional divisor $E_i$ of the blowup of one of the basepoints
  $p_i$ of the net.
\end{proposition} {\bf Proof:} First let us show that
$\operatorname{cont}_R$ must be birational. The key point is that
since $\operatorname{cont}_R$ is the contraction of an extremal ray,
all the curves contracted must be numerical multiples of each
other. 

First suppose $\operatorname{dim} Z =0$. Then all curves on $X$
are numerical multiples of each other, which is clearly false. 

Next suppose $\operatorname{dim} Z =1$. Choose one of the exceptional
divisors $E_i$. I claim that the morphism $\operatorname{cont}_R$
cannot contract any curve in $E_i$. For any such curve is a numerical
multiple of $l_i$, therefore all curves contracted are numerical
multiples of $l_i$, implying that $\operatorname{cont}_R$ is the
contraction of $E_i$, which is birational, contradicting our
hypothesis. So the restriction of $\operatorname{cont}_R$ to $E_i$
contracts no curve, therefore is a surjection $E_i \arrow Z$. But
there can be no surjection from $E_i \iso \P^2$ onto a curve, since
the fibres over distinct points would be disjoint curves in $\P^2$.

Next suppose $\operatorname{dim} Z =2$. As above we get
${(\operatorname{cont}_R)}_{|E_i}: E_i \arrow Z$, a map from $\P^2$ to a
smooth surface which contracts no curves. As before the image cannot
be a curve so it must have dimension $2$: therefore
${(\operatorname{cont}_R)}_{| E_i}$ is surjective. Since
$\operatorname{cont}_R$ contracts no curves, Stein Factorization
\cite[III, Corollary 11.5]{Hartshorne1977} shows that it is a finite
morphism. So the pushforward map ${(\operatorname{cont}_R)}_* :
\Pic(E_i)_\Q \arrow \Pic(Z)_\Q$ is surjective, implying that $\rho(Z)
\leq 1$. On the other hand, \cite[Theorem 3.2]{Mori1982} says that the
contraction of an extremal ray lowers the Picard number by 1, so
$\rho(Z)=\rho(X)-1=8$, which is a contradiction.

So we may assume that $\operatorname{cont}_R$ is birational, and
therefore given by one of the 5 possibilities on Mori's list. My claim
is that only the second of these 5 cases can occur for $X$ the blowup
of $\P^3$ in the base locus of a net of quadrics. To see this, we
will use adjunction for each of the divisors $D$ above. (This is
valid, since each $D$ is a normal divisor in a smooth variety.) For
any curve $C$ contained in $D$, we have $K_X \circ C = (K_X \otimes
\curlyo_D) \circ C$, where the first product is in $CH^*(X)$ and the
second in $CH^*(D)$.  Adjunction then lets us write the second
expression as $(K_D - D) \circ C$. Let us see what this gives in each
of the $5$ cases above, for some choice of $C$.

\begin{enumerate}

\item Fix a section $S_0$ of the bundle $D$, and let $C$ denote any
  fibre. Then one can show \cite[Corollary V.2.11]{Hartshorne1977}
  that $K_D \equiv -2S_0+kC$, for some integer $k$. In particular,
  since $C^2=0$ we have $K_D \circ C = -2$. That gives $K_X \circ C =
  (K_D - D) \circ C = -2- D \circ C$.

  On the other hand, $X$ is the blowing up of a smooth curve in $Z$,
  so we have $K_X= f^*(K_Z)+\curlyo_X(D)$. Since $C$ is contracted by
  $f$, we get $K_X \circ C = D \circ C$. Equating these expressions
  gives $-2-D \circ C= D \circ C$, hence $D \circ C=-1$. Therefore
  $K_X \circ C = -2+1=-1$.  This is impossible in our case, since all
  the coefficients of $K_X$ with respect to the usual basis of
  $\Pic(X)$ are even, hence we must have $K_X \circ C \in 2 \Z$ for
  any curve $C$.

\item This case does occur for our varieties $X$ : blowing down any
  exceptional divisor $E_i$ (where $p_i$ is a basepoint with no
  infinitely near basepoints) gives an example.  I claim these are the
  only examples: any extremal contraction $f: X \arrow Z$ to a smooth
  $Z$ with exceptional divisor $D \iso \P^2$ and $\curlyo_D(D) \iso
  \curlyo(-1)$ must have $D = E_i$ for some basepoint $p_i$. 

  To see this, suppose that $D$ is a divisor satisfying the above
  conditions, distinct from each of the $E_i$. I claim that $D$ must
  be disjoint from each $E_i$. To prove this, suppose $D$ is not
  disjoint from $E_i$: then the intersection $D \cap E_i$ is a curve
  $\Gamma$. Since $f$ contracts $D$, it must contract the curve
  $\Gamma$, hence must contract all of $E_i$, since all curves in
  $E_i$ are numerically equivalent up to constant. Since $D$ is
  irreducible, this gives $D=E_i$.  This contradicts our assumption,
  so we must have $D$ disjoint from $E_i$.  So we can contract all the
  exceptional divisors $E_i$ without changing the isomorphism class of
  $D$ or the normal bundle of $D$. This gives an effective divisor
  $D_0 \subset \P^3$ isomorphic to $\P^2$ with normal bundle
  $\curlyo_{\P^2}(-1)$, which is impossible. So we must have $D=E_i$
  for some basepoint $p_i$.

\item In this case we have $K_D = \curlyo(-2,-2)$, hence $K_D -
  \curlyo_D(D) = \curlyo(-2,-2)-\curlyo(-1,-1)=\curlyo(-1,-1)$. Let
  $C$ be a ruling of $D$: then $K_X \circ C = \curlyo(-1,-1) \circ C =
  -1$. Again this is impossible, since $K_X$ has even coefficients.

\item In this case $\curlyo_D(D) \iso \curlyo_D \otimes \curlyo_{\P^3}
  (-1)$. We can compute $K_D$ using adjunction: viewing $D$ as a
  divisor in $\P^3$, we have $K_D = (K_{\P^3}+ \curlyo_{\P^3}(D))_{|
    D} = (\curlyo(-4) + \curlyo(2))_{| D} = \curlyo(-2)_{| D}$. So
  $K_{X | D} = K_D - \curlyo_D(D) = ( \curlyo (-2)-\curlyo (-1))_{| D}
  = \curlyo(-1)_{| D}$. But then if $C$ is a ruling of the cone, we
  have $K_X \circ C = \curlyo(-1)_{| D} \circ C =-1$.  Again this is
  impossible, since $K_X$ has even coefficients.

\item Here $K_D = \curlyo_{\P^2}(-3)$, and $\curlyo_D(D) =
  \curlyo_{\P^2}(-2)$. Let $C$ be a line in $D \iso \P^2$: then $K_X
  \circ C= (K_D - D) \circ C = \curlyo_{\P^2}(-1) \circ C = -1$. Again
  this is impossible, since $K_X$ has even coefficients.
\end{enumerate}
So as claimed, the only possibility for the contraction of a
$K$-negative extremal ray of $\Curv{X}$ is the contraction of one of
the exceptional divisors $E_i$. \quad QED

We can stretch this argument further:

\begin{corollary}
  Let $X$ be as above, and suppose $X'$ is a smooth SQM of $X$. Then
  the same conclusion holds as for $X$: for any extremal ray $R$ of
  $\Curv{X'}$, the contraction morphism
  $\operatorname{cont}_R$ is of type $2$ on Mori's list.
\end{corollary} {\bf Proof:} Any SQM $\alpha: X' \dashrightarrow X$
with $X'$ smooth induces an isomorphism $N^1(X)_\Z \iso N^1(X')_\Z$
which identifies $K_X$ and $K_{X'}$. In particular $K_{X'}$ is
$2$-divisible in $ N^1(X')_\Z$. The proof of the previous proposition
then applies again.  \quad QED

So suppose we have a smooth SQM $X'$ of $X$, and two $K_{X'}$-negative
extremal rays $R_1$, $R_2$ of $\Curv{X'}$. By the proposition,
the corresponding contractions blow down divisors $D_1$ and $D_2$ in
$X'$, each one a copy of $\P^2$. If $D_1 \cap D_2$ was nonempty, it
would be some curve $C$ say. But then all curves in $D_1$ and all
curves in $D_2$ would be numerical multiples of $C$, which by the cone
theorem implies that the contraction morphism associated to $R_1$ say
must also contract $D_2$. This contradicts the conclusion of Theorem
\ref{theorem-mori3folds} that the exceptional locus is irreducible. We
conclude that any set $D_1, \ldots, D_n$ of such divisors must be
pairwise disjoint. Since each one has normal bundle $\mathcal{O}(-1)$
in $X'$, this implies further that their classes in $N^1(X')$ are
linearly independent. So we can perform a sequence of blowdowns $X' =
X_0 \arrow X_1 \arrow \cdots \arrow X_n$, where $X_i$ is the variety
obtained by contracting $D_1, \ldots, D_i$. The Picard number drops by
$1$ at each stage, and $X_n$ must have Picard number at least $1$, so
we conclude the following:

\begin{corollary} \label{cor-knegextremalrays}
  Suppose $X'$, $D_1, \ldots, D_n$ are as above. Then $n \leq \rho(X')
  - 1$. In particular, $\Curv{X'}$ has at most
  $\rho(X')-1=\rho(X)-1=8$ $K$-negative extremal rays.
\end{corollary}

Now let us restrict to the case where $X'$ is an SQM obtained from $X$
by a sequence of flops. (In fact we will see in the next section that
these are all the SQMs of $X$.) Since $-K_{X'}$ is nef for any $X'$
obtained from $X$ by a sequence of flops, the cone of curves
$\Curv{X'}$ is contained in the closed halfspace $\{ C \in
N_1(X') | K_{X'} \cdot C \leq 0 \}$. So the only other extremal rays
of $\Curv{X'}$ are those in the hyperplane $K^{\perp}$.  The
class of a curve in $X'$ lies in this hyperplane if and only if the
curve is $f'$-fibral, and since the fibres of $f'$ are 1-dimensional,
there are only finitely many classes of such curves. (Indeed if $g:Y
\arrow Z$ is any morphism with 1-dimensional fibres, there are only
finitely many classes of $g$-fibral curves, because the class of a
fibre has only finitely many decompositions in the monoid of effective
classes in $N_1(Y)_\Z$.) All this strongly suggests that each of the
nef cones $\Nef{X'}$ should be rational polyhedral. However, it is
{\it a priori} possible that the cone behaves strangely in a
neighbourhood of $K^{\perp}$ yielding extremal rays which are not
spanned by the class of any curve.

We prove that this bad behaviour does not occur for $X$: in other
words, that $X$ has rational polyhedral cone of curves. Under the
additional assumption that the Mordell--Weil rank $\rho$ is $7$ (or
equivalently that the net has no reducible member) we prove the same
conclusion for any $X'$ obtained from $X$ by flopping a set of fibral
curves. 

\begin{theorem} \label{theorem-nefcones}
Suppose $X$ is the blowup of $\P^3$ in the base locus of a net of
quadrics with 8 distinct basepoints. Then the closed cone of curves
$\Curv{X}$ is rational polyhedral, spanned by the classes
$l_i$ of lines in the exceptional divisors $E_i$ together with the 28
classes $C_{ij}$. If in addition the net has Mordell--Weil rank
$\rho=7$ then $\Curv{X'}$ is rational polyhedral for any
$X'$ obtained from $X$ by flopping a set of fibral curves.

Dually, the nef cones $\Nef{X}$ and $\Nef{X'}$ are rational polyhedral
in the situations described.
 
\end{theorem}

{\bf Proof:} Corollary \ref{cor-knegextremalrays} showed that
$\Curv{X}$ has only finitely many $K$-negative extremal
rays, so it suffices to show there are only finitely many extremal
rays in $K^\perp$.

Consider a divisor class of the form $D_{ij}=H-E_i-E_j$
on $X$. This class is represented by the proper transform on $X$ of
any plane in $\P^3$ passing through the points $p_i$ and $p_j$ so its
base locus is the curve $C_{ij}$. Therefore if $C$ is any irreducible
curve on $X$ which is not one of the curves $C_{ij}$, we must have
$D_{ij} \cdot C \geq 0$. So all but finitely many irreducible curves
$C$ on $X$ satisfy $D_{ij} \cdot C \geq 0$ for all $i, \ j$. In
particular, any limit ray $R$ of a sequence of irreducible curves
which is not contained in the cone spanned by the $C_{ij}$ must
satisfy $D_{ij} \cdot R \geq 0$ for all $i, \ j $. We
know all the extremal rays of $\Curv{X}$ except those in
$K^\perp$, so any other extremal ray $R$ must also satisfy $K_X \cdot
R=0$. By computation the cone defined by the inequalities $D_{ij}
\cdot C \geq 0$ and $K_x \cdot C =0$ is spanned by a finite set of
vectors of the form $nl-(n-1)l_{i_1}-l_{i_2}-\cdots-l_{i_{n+2}}$, for
$n=2,\ldots,6$. Therefore this cone is contained in the cone $\R_+ \{
C_{ij} \}$ spanned by the classes $C_{ij}$. This proves that
$\Curv{X} \cap K^\perp = \R_+ \{C_{ij} \}$, and therefore
$\Curv{X}$ is a rational polyhedral cone whose extremal rays
are spanned by the classes $l_i$ and $C_{ij}$. This proves the first claim.

Now suppose $\rho=7$, and let us prove the claim about the cones
$\Curv{X'}$. The idea is similar to the proof of the claim
about $\Curv{X}$, but the role of the divisors $D_{ij}$ is
now played by an infinite set of movable divisors. Again Corollary
\ref{cor-knegextremalrays} tells us that $\Curv{X'}$ has
only finitely many $K$-negative extremal rays, so it suffices to prove
that there are only finitely many extremal rays in $K^\perp$.

Recall that in general the action of $\Pic^0(X_\eta)$ on $N^1(X)$ is
given by the formula $\psi_y(x) = x + (x \cdot F)y + V(x,y)$, where $x
\in N^1(X)$, $y \in \Pic^0(X_\eta)$, and $V(x,y)$ is a vertical
divisor. In the case $\rho=7$, the only vertical divisors are
multiples of $-\frac12K_X$, so we get $\psi_y(x) = x + (x \cdot F)y
+m(-\frac12K_X)$. In particular if $x=D_{ij}$ and $y=n(E_k-E_l)$ we
have $D_{ij} \cdot F= 2$ and hence $\psi_y(D_{ij})= H-E_i-E_j
+2n(E_k-E_l) + m(-\frac12K_X)$. Also, since the base locus of $D_{ij}$
is the curve $C_{ij}$ and $\Pic^0(X_\eta)$ acts by
pseudo-automorphisms over $\P^2$, the base locus of any such divisor
$\psi_y(D_{ij})$ is a finite union of fibral curves. If we then flop
some fibral curves to obtain $X'$, the base locus of the proper
transform $\psi_y(D_{ij})'$ is again a finite union of fibral
curves. The upshot is that for any irreducible curve $C$ on $X'$,
either $\psi_y(D_{ij})' \cdot C \geq 0$ for all $i,\, j$ and all $y$
in $\Pic^0(X_\eta)$, or else $C$ is one of the finitely many classes
of fibral curves on $X'$.

Now suppose $R$ is an extremal ray of $\Curv{X'}$ which lies
in the subspace $K^\perp$. As before, any limit ray $R$ of a sequence
of irreducible curves which is not in the cone spanned by the classes
of fibral curves must satisfy $\psi_y(D_{ij})' \cdot R \geq 0$ for all
$i,\, j$ and all $y$ in $\Pic^0(X_\eta)$. Suppose that $R$ is such a
ray. Since $R \subset K^\perp$, any class $C$ which spans $R$ has $K
\cdot C =0$, implying $(\psi_y(D_{ij}))' \cdot C = (D_{ij} + 2y)'
\cdot C \geq 0$ for any $x \in N^1(X)$ and $y \in \Pic^0(X_\eta)$. In
particular if we put $x=D_{ij}$ and $y=n(E_k-E_l)$ we get
$((H-E_i-E_j) +2n(E_k-E_l))' \cdot C \geq 0$ for all indices
$i$, $j$, $k$, $l$ and all integers $n$. Now if $C=al+\sum_ib_il_i$ with
coefficients $b_i$ not all equal, then $(E_k-E_l)' \cdot C <0$ for
some $k$ and $l$, implying $((H-E_i-E_j) +2n(E_k-E_l))' \cdot C < 0$
for some indices $i$, $j$, $k$ ,$l$ and $n$ sufficiently large. This
contradicts our choice of $R$. So the only possibility is that all
coefficients $b_i$ are equal, which implies that $R$ is the ray
spanned by $4l-\sum_il_i$, the class of a fibre. We conclude that
$\Curv{X'} \cap K^\perp$ is spanned by the classes of fibral
curves, which are finite in number, and therefore that
$\Curv{X'}$ is a rational polyhedral cone, as claimed. \quad QED

The cone conjecture concerns the nef effective cone $\Nef{X}^e$ rather
than the whole nef cone. However, in our situation these cones
coincide:

\begin{proposition} \label{prop_effective} Suppose $X$ is the blowup
  of $\P^3$ in the base locus of a net of quadrics with 8 distinct
  basepoints, and $X'$ is obtained from $X$ by a sequence of flops of
  fibral curves. Then any nef divisor on any of the varieties $X'$ is
  semi-ample, hence effective. In other words, $\Nef{X'}^e=\Nef{X'}$
  for all such $X'$.
\end{proposition}

{\bf Proof:} Since $\Nef{X'}$ is rational polyhedral it suffices to
prove that any integral divisor in $\Nef{X'}$ is effective. The
proposition is immediate for multiples of $-\frac12K_{X'}$ so assume
$D$ is a nef integral divisor which is not such a multiple. By the
Basepoint-Free theorem \cite[Theorem 3.3]{KollarMori1998} it suffices
to show that the divisor $D-\frac12K_{X'}$ is big, which by the
numerical criterion for bigness of nef divisors \cite[Theorem
  2.2.16]{Lazarsfeld2004} is equivalent to
$(D-\frac12K_{X'})^3>0$. Now $(-\frac12 K_{X'})^3=0$ and $D^3 \geq 0$
since $D$ is nef; also $D^2 \cdot (-\frac12K_{X'}) \geq 0$ since $D^2
\in \Curv{X'}$. So it suffices to prove that $D \cdot
(-\frac12 K_{X'})^2>0$. But $(-\frac12 K_{X'})^2$ is the class $F$ in
$N^1(X')$ of any fibre of $f'$. If $D \cdot F=0$ then since $D$ is nef
we must have $D \cdot C =0$ for $C$ the class of any $f'$-fibral
curve. The classes of such curves span the codimension-1 subspace
$K^\perp$ of $N_1(X')$ so $D$ is a multiple of $-\frac12K_{X'}$,
contradicting our initial assumption. Therefore $D \cdot (-\frac12
K_{X'})^2 >0$ as required. \quad QED

We have proved the first statement of Theorem \ref{theorem-mainthm},
namely that $\Nef{X}$ is a rational polyhedral cone spanned by
effective divisors. However, the first prediction of Conjecture
\ref{conj-totaro} does not seem to follow immediately. The conjecture
predicted there should be a rational polyhedral fundamental domain for
the action of $\PsAut^*(X,\Delta)$ on $\Nef{X}$. (Recall that $\Delta$
is a $\Q$-divisor $\frac12 D$ for some smooth member $D$ of the linear
system $|-2K_X|$.) To verify that statement for $X$, we use the
following theorem of Looijenga \cite[Proposition 4.1, Application
  4.15]{Looijenga2009}. (We state a stronger form than we need at
present, for use in the next section.)

\begin{theorem}[Looijenga] \label{theorem-looijenga}
Let $V$ be a real vector space with $\Z$-structure and $C$ a strictly
convex open cone in $V$ with nonempty interior. Let $G$ be a subgroup
of $GL(V_\Z)$ which preserves $C$. Suppose there is a rational
polyhedral cone $U$ in $\overline{C}$ such that $G \cdot U$
contains $C$. Then $G \cdot U$ is equal to the convex hull $C_+$ of
the rational points in $\overline{C}$, and there exists a rational
polyhedral fundamental domain for the action of $G$ on $C_+$.
\end{theorem}

\begin{corollary}
The first statement of Conjecture \ref{conj-totaro} holds for $X$:
there is a rational polyhedral fundamental domain for the action of
$\Aut^*(X,\Delta)$ on $\Nef{X}^e$.
\end{corollary} {\bf Proof:} We have just seen that
$\Nef{X}^e=\Nef{X}$, a rational polyhedral cone. Applying Theorem
\ref{theorem-looijenga} with $\overline{C}=U=\Nef{X}$ and
$G=\Aut^*(X,\Delta)$ we get the result.  \quad QED

\begin{corollary}
For $X$ as above, the group $\Aut^*(X)$ is finite.
\end{corollary} {\bf Proof:} The cone $\Nef{X}$ is a
  rational polyhedral cone preserved by $\Aut^*(X)$. I claim that an
  infinite subgroup $G$ of $GL(N^1(X)_\Z)$ cannot preserve a strictly
  convex rational polyhedral cone with nonempty interior. For the
  action of any element of $g \in G$ must permute the primitive
  integral vectors in the extremal rays of the cone, and this
  permutation determines the action of $g$. So $G$ is realised as a
  subgroup of a finite permutation group. \quad QED

\section{Movable cone} \label{section-movablecone}

The aim of this section is to prove the second part of Theorem
\ref{theorem-mainthm}: if $X$ is the blowup of $\P^3$ in the base
locus of a net of quadrics with 8 distinct basepoints, and the
Mordell--Weil rank $\rho(X)$ equals 7, then there is a rational
polyhedral fundamental domain for the action of $\PsAut^*(X)$ on
$\Mov{X}^e$. 

We remark that the second part of Theorem \ref{theorem-mainthm} is
more difficult than the first: in particular, the effective movable
cone $\Mov{X}^e$ is in general not rational polyhedral in the context
we are considering. To see this, recall from the introduction that the
Mordell--Weil group $\Pic^0(X_\eta)$ of the generic fibre of $f$ acts
on $X$ by pseudo-automorphisms. I claim that the representation
$\Pic^0 (X_\eta) \arrow GL(N^1(X)_\Z)$ is faithful. To see this, note
that any rational section $D$ of $f$ has an open subset covered by
rational curves $C$ with $D\cdot C=-1$ (images under the section of
lines in $\P^2$). If $D_0$ and $D_1$ are different rational sections,
then the element $D_1-D_0 \in \Pic^0 (X_\eta)$ maps $D_0$ to $D_1+V$
where $V$ is an effective divisor pulled back from $\P^2$. Then for
all but finitely many curves $C \subset D_0$ we have $(D_1+V)\cdot C
\geq 0$, so $D_1+V$ and $D_0$ are numerically distinct, and therefore
$D_1-D_0$ is not in the kernel of the representation. This proves the
claim, and we conclude that if the Mordell-Weil group is infinite,
then so too is its image in $GL(N^1(X)_\Z)$. Pseudo-automorphisms of
$X$ preserve the effective movable cone $\Mov{X}^e$, and (as mentioned
before) an infinite subgroup of $GL(N^1(X)_\Z)$ cannot act on a
strictly convex rational polyhedral cone with nonempty
interior. Therefore the effective movable cone cannot be rational
polyhedral unless the Mordell--Weil group is finite (or equivalently,
as explained in Section \ref{sect-prelims}, the net contains the
maximum number 7 of reducible members).

The structure of the proof is as follows. First we show in Proposition
\ref{prop-movablecone} that the cone $\Mov{X}^e$ decomposes as the
union of nef effective cones of SQMs of $X$ which are obtained by
flopping curves in the fibres of $f: X \arrow \P^2$, and the interiors
of these nef cones are disjoint. Moreover (Lemma
\ref{lemma-permutation}), pseudo-automorphisms of $X$ act by permuting
the cones. We saw in the previous section (under the assumption of
maximum Mordell--Weil rank) that each of the nef cones is rational
polyhedral, so it seems reasonable that some finite union of these
cones might provide the fundamental domain we seek. Precisely, by
Theorem \ref{theorem-looijenga}, it is enough to show that the
translates by pseudo-automorphisms of a finite union of these cones
covers the effective movable cone.

To prove this, we again use the elliptic fibration structure on
$X$. As previously mentioned, the Mordell--Weil group $\Pic^0(X_\eta)$
of the generic fibre of the fibration is a subgroup of the
pseudo-automorphism group of $X$. We study the action of this subgroup
on the quotient space $N^1(X/\P^2)$. By general results of Kawamata on
3-dimensional elliptic fibrations (Lemma \ref{lemma-kawamata}), this
action is easy to understand explicitly. Also, in Lemma
\ref{lemma-relativemovable} we are able to compute `by hand' the
relative movable cone, using the explicit geometry of the fibration.

Putting these facts together we find in Lemma \ref{lemma_funddomain} a
rational polyhedral cone in the relative movable cone whose
$\Pic^0(X_\eta)$-translates cover the whole cone. The key point of our
method is to lift this action to the absolute movable cone. In Theorem
\ref{theorem-funddomain} to find a rational polyhedral cone in
$N^1(X)$ whose $\Pic^0(X_\eta)$-translates cover the whole effective
movable cone. Since the Mordell--Weil group is a subgroup of the
pseudo-automorphism group, the $\PsAut(X)$-translates of that rational
polyhedral cone also cover the whole effective cone, as required.

As a final remark before starting the proof, we re-emphasise that most
of our proof is valid for nets of arbitrary Mordell--Weil rank. We
only impose the restriction to nets of maximum rank starting from
Lemma \ref{lemma-rank7results}. It seems likely that the remaining
details can be filled in to give a complete proof for nets of
arbitrary rank. The proof below would be simplified if we restricted
to nets of maximum rank from the start; however, we retain the general
argument as far as possible, to illustrate the fact that our methods
still give a good deal of information in the general case.

Now let us begin the proof. The first step is to show that any movable
divisor on $X$ can be made nef by a sequence of flops. For $D$ a
$\Q$-Cartier $\Q$-divisor on a normal projective variety $Y$, define a
{\it $D$-flopping contraction} of $Y$ to be a proper birational
morphism $f:Y \arrow Z$ to a normal variety $Z$ such that the
exceptional set of $f$ has codimension at least 2 in $Y$, the
canonical class $K_Y$ is numerically $f$-trivial, and $-D$ is
$f$-ample. The {\it $D$-flop} of $f$ is then defined to be the
$(K_Y+D)$-flip of $f$. We need the following result \cite[Theorem
  6.14, Corollary 6.19]{KollarMori1998}:

\begin{lemma} \label{lemma-flops}
Suppose $Y$ is a threefold with terminal singularities and $f: Y
\arrow Z$ a $D$-flopping contraction (for some $\Q$-divisor $D$). Then
the $D$-flop of $f$ exists. Moreover, any sequence of extremal
$D$-flops on a terminal threefold is finite.
\end{lemma}

Here an {\it extremal} flop is one for which the flopping contraction
has relative Picard number 1. In particular if $(X,\Delta)$ is a
$\Q$-factorial klt pair with $\Delta$ effective and $f$ is the
contraction of a $(K_X+\Delta)$-negative extremal ray, it is extremal,
because all the curves contracted are numerical multiples of each
other, by the cone theorem.

This lemma enables us to find the desired decomposition of the movable
cone of $X$ into nef cones:
\begin{proposition} \label{prop-movablecone} The effective movable cone
  $\Mov{X}^e$ decomposes as a union of nef cones of SQMs of $X$:
\begin{align*}
  \Mov{X}^e = \bigcup  \Nef{X',\alpha} 
\end{align*}
where the union on the right hand side is over all SQMs $\alpha:X'
\dashrightarrow X$. All these SQMs are obtained by flopping fibral
curves. The interiors of the cones $\Nef{X',\alpha}$ are disjoint.
\end{proposition}
{\bf Proof:} Suppose $D \in \Mov{X}^e$ is an effective $\Q$-divisor on
$X$ which is not nef. We know a movable divisor cannot be negative on
a $K$-negative extremal ray, so the same is true of a divisor in
$\Mov{X}$.  So the description of $\Curv{X}$ tells us that $D
\cdot C_{ij} <0$ for some $i, \, j$. Choosing some $\epsilon >0$
sufficiently small, the cone theorem for the klt pair $(X, \epsilon
D)$ tells us that the contraction of $C_{ij}$ exists. By Lemma
\ref{lemma-flops} we can perform a flop to obtain a variety $X'$ on
which $D' \cdot C' >0$, where $D'$ is the proper transform of $D$, and
$C'$ the cocentre of the flop. If $D'$ is now nef we can stop. If not
then $D'$ is negative on some extremal ray of
$\Curv{X'}$. Again $D'$ cannot be negative on a $K$-negative
extremal ray so $D' \cdot R <0$ for some extremal ray $R \in
K^{\perp}$. Choosing $\epsilon'>0$ sufficiently small again the cone
theorem for $(X',\epsilon'D')$ tells us that $R$ is spanned by the
class of a curve $C'$ (necessarily a component of a fibre) and the
$D'$-flop exists. Continuing in this way, the lemma above guarantees
that this process terminates at some finite stage: that is, the proper
transform of $D$ becomes nef after some finite sequences of flops of
fibral curves.

If $D$ is not a $\Q$-divisor, the same argument works
by adding at each stage a sufficiently small ample $\R$-divisor $D_1$
so that $D'+D_1$ is a $\Q$-divisor which has negative degree on the
same fibral curves as $D'$.

We have shown that any effective movable divisor belongs to one of the
effective nef cones $\Nef{X'}^e$ where $X'$ is obtained from $X$ by
flopping fibral curves. So we have the inclusion $\Mov{X}^e \subset
\cup \Nef{X',\alpha}^e$. The reverse inclusion is clear, since an ample
divisor on $X'$ is movable on $X$, so taking closures and intersecting
with the effective cone we get $\cup \Nef{X',\alpha}^e \subset \Mov{X}^e
$. Finally by Proposition \ref{prop_effective} we get the statement
above. 

To see that these flops give all the SQMs of $X$ up to isomorphism, suppose
that $\beta: Y \dashrightarrow X$ is any SQM. By the argument above we
have $\Nef{Y,\beta} \subset \cup \Nef{X',\alpha}$, therefore the ample
cone of $Y$ must intersect the ample cone of one of the flops, say
$\alpha_i: X_i \dashrightarrow X$. So there exists a divisor $D$ on
$X$ such that ${\alpha_i}^{-1}_* D$ and $\beta^{-1}_* D$ are ample on
$X_i$ and $Y$ respectively. Therefore

\begin{align*}
X_i = \operatorname{Proj} R(X_i,{\alpha_i}^{-1}_* D) \iso
\operatorname{Proj} R(Y,\beta^{-1}_*
D) = Y
\end{align*}
and the isomorphism is compatible with $\alpha_i$ and $\beta$. 

Finally, the same argument applied to $2$ SQMs $\alpha_1: X_1
\dashrightarrow X$ and $\alpha_2: X_2 \dashrightarrow X$ shows that
the interiors of the cones $\Nef{X_1,\alpha_1}$ and
$\Nef{X_2,\alpha_2}$ are disjoint in $\Mov{X}^e$. \quad QED

This decomposition of the movable cone is compatible with the action
of pseudo-automorphisms:
\begin{lemma} \label{lemma-permutation} The group $\PsAut(X)$ acts on
  $\Mov{X}^e$ by permuting the nef cones of the small modifications of
  $X$. More precisely suppose $\phi \in \PsAut(X)$
  and $\alpha:Y \dashrightarrow X$ is an SQM of $X$. Then
  $\phi_*(\Nef{Y})=\Nef{Y'}$ for $\alpha':Y' \dashrightarrow X$
  some other SQM of $X$.
\end{lemma} {\bf Proof:} Suppose $\alpha: Y \dashrightarrow X$ is the
marking of $Y$. Then $\alpha'$ is the SQM given by $\phi \circ \alpha: Y
\dashrightarrow X$. To see this works take any divisor $D$ on $X$ such
that $\alpha^{-1}_*(D)$ is nef on $Y$: in other words $D$ belongs to
$\Nef{Y,\alpha}$.  Then putting $\Delta=\phi_*(D)$ we have
$(\phi \circ \alpha)^{-1}_*(\Delta)=\alpha^{-1}_*(D)$ nef on $Y$. So $\Delta$
belongs to $\Nef{Y,\alpha'}$ and $\phi_*$ maps $D$ to $\Delta$.
This holds for any $D \in \Nef{Y}$ and so
$\phi_*(\Nef{Y,\alpha}) \subset \Nef{Y,\alpha'}$. Exchanging
$\phi$ and $\phi^{-1}$ we get the result. \quad QED

The next step in the proof is to study the action of
pseudo-automorphisms (more precisely, the Mordell--Weil group) on the
relative movable cone $\Mov{X/\P^2}^e$. First we must understand the
relationship between the vector spaces $N^1(X)$ and
$N^1(X/\P^2)$. Recall that $N^1(X/\P^2)$ was defined as the space of
Cartier divisors on $X$ with real coefficients modulo numerical
equivalence over $\P^2$, where a divisor $D$ is numerically trivial
over $\P^2$ if $D \cdot C =0$ for every curve $C$ which maps to a
point on $\P^2$. Those curves span a subspace of $N_1(X)$, so dually
there is a projection map $p$ from $N^1(X)$ to $N^1(X/\P^2)$. As a
piece of notation, from now on we will write $[D]$ to denote the image
$p(D) \in N^1(X/\P^2)$ of an element $D \in N^1(X)$.

The first important question is to characterise this projection map:
\begin{lemma} \label{lemma-projection}
The projection $p: N^1(X) \arrow N^1(X/\P^2)$ has
$1$-dimensional kernel spanned by the class of $-\frac{1}{2}K_X$.
\end{lemma} {\bf Proof:} A class $D \in N^1(X)$ maps to $0$ in
$N^1(X/\P^2)$ if and only if $D \cdot C=0$ for every fibral curve $C$
on $X$. The description of $\Curv{X}$ in the previous section shows
that the classes of fibral curves span the hyperplane $K^\perp$, so
$D$ must be a multiple of $-\frac12K_X$. \quad QED

For later use we introduce some notation related to this projection:
for a class $D \in N^1(X)$, we denote its image $p(D) \in N^1(X/\P^2)$
by $[D]$.

The following lemma of Kawamata \cite[Lemma 3.5]{Kawamata1997} shows
that the action of $\Pic^0(X_\eta)$ is easy to understand if we pass
to a suitable quotient space of $N^1(X/\P^2)$.

\begin{lemma} \label{lemma-kawamata} Let $V(X/\P^2)$ denote the
  subspace of vertical divisors of $f:X \arrow \P^2$ and $W(X/\P^2)$
  the affine subquotient space $\left\{ x \in N^1(X/\P^2)/V(X/\P^2) : x
  \cdot F = 1 \right\}$. Then $\Pic^0(X_\eta)$ acts properly
  discontinuously on $W(X/S)$ as a group of translations, and has
  fundamental domain a rational polyhedron $\Pi$.
\end{lemma}

Next we compute the relative effective and movable cones. It turns out
that the fibration $f$ determines these cones in the simplest way one
could hope for. That is, certain curves in the fibres of $f$ give
obvious classes in $N_1(X/\P^2)$ on which any $f$-effective or
$f$-movable divisor must have nonnegative degree, and it turns out
that these obvious classes actually suffice to determine the cones
completely. The precise result is the following:

\begin{lemma} \label{lemma-relativemovable} Suppose $f:X \arrow \P^2$
  is as before. For $i=1, \, 2$ let $D^a_i$ (resp. $F^a_i$) denote the
  components of the reducible quadric $Q_i$ in the net
  (resp. components of the fibre $f^{-1}(\eta_i)$ where $\eta_i$ is
  the generic point of $f(Q_i)$) and let $F$ denote the class of a
  fibre of $f$.  Then
\begin{align*}
  B^e(X/\P^2) &=\{x \in N^1(X/\P^2) : x \cdot F >0 \} \cup
  \R_+\{[D^a_i]\} \cup \{0\} \\ \Mov{X/\P^2}^e &= \{x \in N^1(X/\P^2)
  : x \cdot F >0, \, x \cdot F^a_i \geq 0 \text{ for all } a, \,
  i\} \cup \{0\}.
\end{align*}
Note that two divisors whose classes are equal in $N^1(X/\P^2)$ differ
by a multiple of $-\frac12K_X$, so intersection numbers with all $f$-fibral
curves are well-defined.
\end{lemma} {\bf Proof:} First suppose that $D$
is an $f$-effective Cartier divisor with degree $k \leq 0$ on the
generic fibre.  There exists a nonempty open set $U \subset \P^2$ such
that $D(f^{-1}(U)) \neq 0$. Choose a nonzero section $s \in
D(f^{-1}(U))$. Then the class of the divisor
$\Delta=\overline{\{s=0\}}$ differs from $D$ only on the codimension-1
subset $X \backslash f^{-1}(U)$. In particular these classes have the
same degree $k$ on the generic fibre. Since $\Delta$ is effective this
implies $k=0$.  Moreover $k=0$ implies that $\Delta$ is a sum of
vertical divisors so its class in $N^1(X)$ belongs to the cone $V$
spanned by $-\frac12 K_X$ and the $D^a_i$. Finally the divisor
$\Delta-D$ is supported in $X \backslash f^{-1}(U)$ therefore its
support maps onto a curve in $\P^2$. So any divisor in the support of
$\Delta-D$ must also have class in the cone $V$. So for any
$f$-effective class $D$ with degree $\leq 0$ on the generic fibre we
can write $D=V_1-V_2$ where $V_i$ are classes in $V$. The image of $V$
in $N^1(X/\P^2)$ is the cone $\R_+\{[D^a_i]\}$, which is closed under
negation since for any $i$ we have $[D^1_i]+[D^2_i]=0$ in
$N^1(X/\P^2)$. Therefore $[D]=[V_1]-[V_2]$ belongs to this cone as
claimed. This proves that the left-hand side of the first equation is
contained in the right-hand side. To prove the reverse inclusion,
first note that if a divisor $D$ has positive degree on an irreducible
fibre $F$ then the restriction $D_{|F}$ is ample hence effective. But
standard results on semicontinuity of cohomology \cite[Corollary
  III.12.9]{Hartshorne1977} show that any section of $D_{|F}$ is the
restriction of a section in $D(f^{-1}(U))$ for $U \in \P^2$ some open
subset. By definition that means $D$ is $f$-effective. Finally, all
divisors in the cone $V$ are effective by definition hence
$f$-effective, so all elements of $\R_+\{[D^a_i]\}$ lie in the
$f$-effective cone. This completes the proof of the claim about the
$f$-effective cone.

Now we must prove the claim about the $f$-movable cone. First note
that if $D$ is a Cartier divisor in $N^1(X)$ with $D \cdot F^a_i <0$
for some $a$ and $i$, then $D$ cannot be $f$-movable. For suppose $C$
is a curve in $X$ with class $F^a_i$. If there was an open set $U
\subset \P^2$ containing the point $f(C)$ and a section of
$D(f^{-1}(U))$ not vanishing identically along $C$ we would have $D
\cdot C \geq 0$ contradicting our assumption: therefore every such
curve $C$ is contained in $\operatorname{Supp \ Coker } (f^*f_*
\mathcal{O}_X(D) \arrow \mathcal{O}_X(D))$. Since these curves $C$
fill up the divisor $D^a_i$ we conclude that $D$ cannot be
$f$-movable. So the $f$-movable cone is contained in the cone $\{x
\cdot F^a_i \geq 0 \text{ for all } a, \, i \}$. If moreover $D$ is a
nonzero $f$-effective $f$-movable class, it must have $D \cdot
F>0$. For otherwise by the description of $B^e(X/\P^2)$ we would have
$[D] \in \R_+\{[D^a_i]\}$. Any nonzero point in this cone has the form
$\sum r_{ai} [D^a_i]$ with $r_{1i}$ and $r_{2i}$ not equal for all
$i$. Say $r_{1i} > r_{2i}$: then $[D] \cdot F^1_i <0$, contradicting
our previous conclusion. So we have shown that left-hand side is
contained in the right-hand side in the second equality above.

Conversely suppose that $x \in N^1(X/\P^2)$ satisfies $x \cdot
F^a_i \geq 0$ for all $a,\, i$ and $x \cdot F>0$: we want to show
that $x$ belongs to the $f$-effective $f$-movable cone. First note
that any such $x$ is $f$-effective by our description of the
$f$-effective cone. Next suppose that $D$ is a divisor class with $D
\cdot F^a_i >0$ for each $i$. Since $F=F_i^1+F_i^2$ for any $i$, the
restriction of such a $D$ to any irreducible fibre is ample. Also
since $D \cdot F^a_i>0$ for each $a$ and $i$ the restriction of $D$ to
any component of a fibre not containing a line $C_{ij}$ is ample. So
by taking a sufficiently large multiple $mD$ we get a line bundle
whose restriction to all but finitely many fibres is very ample hence
basepoint-free. Again by the semicontinuity result mentioned above any
section of a line bundle $D_{|F}$ comes from a section in
$D(f^{-1}(U))$ for some open $U \subset \P^2$ containing $f(F)$. So
$\operatorname{Supp \ Coker } (f^*f_* \mathcal{O}_X(D) \arrow
\mathcal{O}_X(D))$ does not contain any point in any of these fibres
and therefore has codimension at least $2$. Therefore the class $[D]$
belongs to the $f$-movable cone. To complete the proof we observe that
any class in the cone $\{x \in N^1(X/\P^2) : x \cdot F >0, \, x
\cdot F^a_i \geq 0 \text{ for all } a, \, i\}$ is the limit of classes
$[D_\alpha]$ with $[D_\alpha] \cdot F^a_i >0$ for all $a$ and $i$. We
have just proved that each class $[D_\alpha]$ belongs to the closed
cone $\Mov{X/\P^2}$ and therefore so does their limit $x$. \quad QED

We use this description of $\Mov{X/\P^2}^e$ together with Lemma
\ref{lemma-kawamata} to find a rational polyhedral cone whose
translates by the Mordell--Weil group cover the relative movable cone:

\begin{lemma} \label{lemma_funddomain}
There is a rational polyhedral subcone $K$ of $\Mov{X/\P^2}^e$ such
that $\Pic^0(X_\eta) \cdot K = \Mov{X/\P^2}^e$. 
\end{lemma} {\bf Proof:} Let $W'(X/\P^2)$ denote
the affine subspace $\left\{y \in N^1(X/\P^2) : y \cdot F = 1 \right\}$
and denote by $q$ the quotient map $W'(X/\P^2) \arrow
W(X/\P^2)$. By definition of the quotient action of $\Pic^0(X_\eta)$,
for any $\phi \in \Pic^0(X_\eta)$ and $x \in N^1(X/\P^2)$ we have
$\phi(q(x))=q(\phi(x))$. By Lemma \ref{lemma-kawamata} the action of
$\Pic^0(X/\P^2)$ on $W(X/\P^2)$ has fundamental domain a rational
polyhedron $\Pi$, and hence for the action on $W'(X/\P^2)$ we have
$\Pic^0(X/\P^2) \cdot q^{-1}(\Pi) = W(X/\P^2)$. Since the
action of $\Pic^0(X/\P^2)$ preserves the $f$-effective $f$-movable
cone, we can intersect with that cone on both sides to get
$\Pic^0(X/\P^2) \cdot (q^{-1}(\Pi) \cap \Mov{X/\P^2}^e) =
\Mov{X/\P^2}^e \cap W(X/\P^2)$. Finally since $\Pic^0(X/\P^2)$ acts
linearly we can multiply on both sides by positive scalars to get
$\Pic^0(X/\P^2) \cdot \R_+ (q^{-1}(\Pi) \cap \Mov{X/\P^2}^e) =
\Mov{X/\P^2}^e$. So taking $K =\R_+ (q^{-1}(\Pi) \cap
\Mov{X/\P^2}^e)$ it remains to show that $q^{-1}(\Pi) \cap
\Mov{X/\P^2}^e$ is a rational polyhedron in $W'(X/\P^2)$. Since $\Pi$
is a rational polyhedron and by Lemma \ref{lemma-relativemovable} the
cone $\Mov{X/\P^2}^e$ is defined by a finite set of inequalities, we
need to show that $q^{-1}(\Pi) \cap \Mov{X/\P^2}^e$ is
bounded. Choosing a section $s$ of $q$ we can write
$W'(X/\P^2)=V(X/\P^2) + \operatorname{im} s$. Let $\Pi'$ denote
the polyhedron $s(\Pi)$: then $q^{-1}(\Pi) = V(X/\P^2) + \Pi'
\subset W'(X/\P^2)$. So suppose a vector $v+s$ with $v \in V(X/\P^2)$
and $s \in \Pi'$ belongs to $q^{-1}(\Pi) \cap \Mov{X/\P^2}^e$. By
Lemma \ref{lemma-relativemovable} the intersection numbers $(v+s)
\cdot F^a_i$ must be nonnegative for all $i$ and $j$. Now $s \cdot
F^a_i$ is bounded for $s \in \Pi'$ by compactness of $\Pi'$, so $v
\cdot F^a_i$ is bounded below for all $a$ and $i$. But $v \cdot F^1_j
= -v \cdot F^2_j$ for any $v \in V(X/S)$, therefore $v \cdot F^a_i$ is
bounded above and below. If we write $v=\sum a^k_l D^k_l$ we have $v
\cdot F^a_i = -2a^i_j$ so the coefficients of $v$ are bounded. So the
subset $q^{-1}(\Pi) \cap \Mov{X/\P^2}^e$ is bounded, hence rational
polyhedral, as required. \quad QED

So far we have considered the action of the Mordell--Weil group on the
relative movable cone. To complete the proof, we lift our result to
the absolute movable cone $\Mov{X}^e \subset N^1(X)$.  Precisely, we
use the previous lemma to find a polyhedral cone in $\Mov{X}^e$ whose
translates by the Mordell--Weil group cover the whole cone. Applying
Theorem \ref{theorem-looijenga} we will see this implies the existence
of a rational polyhedral fundamental domain, thereby completing the
proof of the second part of Theorem \ref{theorem-mainthm}.

We note that so far everything we have said in this section applies to
nets of arbitrary Mordell--Weil rank. Only now do we restrict to the
case of Mordell--Weil rank $\rho=7$.

\begin{lemma} \label{lemma-rank7results}
 Suppose $X$ is the blowup of $\P^3$ in the base locus of a net of
 quadrics which has 8 distinct basepoints, with Mordell--Weil rank
 $\rho=7$. Then $\Mov{X/\P^2}^e$ is the open half-space $\{[D] \in
 N^1(X/\P^2) : [D] \cdot F >0 \}$. Moreover, let $\Pi$ denote the
 rational polyhedron
\begin{align*}
\Pi = \{ \sum_i \alpha_i [E_i] : 0 \leq \alpha_i \leq 1 \text{ for }
i=2,\ldots,8 \, ; \sum_i \alpha_i =1 \}.
\end{align*}
contained in the affine hyperplane $W'(X/\P^2)= \{[D] \in N^1(X/\P^2)
: [D] \cdot F =1 \}$. Then $K=\R_+ \Pi$ is a rational polyhedral cone
satisfying the conclusion of Lemma \ref{lemma_funddomain}: that is,
$\Pic^0(X_\eta) \cdot K = \Mov{X/\P^2}^e$.

\end{lemma} {\bf Proof:} The first claim follows directly from Lemma
\ref{lemma-relativemovable}, since if the Mordell--Weil rank equals 7 there
are no vertical divisors on $X$ other than multiples of $-\frac12K_X$.

For the second claim Lemma \ref{lemma-kawamata} tells us that
$\Pic^0(X_\eta)$ acts on the affine hyperplane $W'(X/\P^2)$ as a group
of translations. Now $\Pic^0(X_\eta)$ has a subgroup $G$ of finite
index generated by the elements $E_j-E_1 \ (j=2,\ldots,8)$. It is
clear that $\Pi$ is a fundamental domain for the action of $G$ on
$W'(X/\P^2)$. Furthermore, $W'(X/\P^2)$ generates $\Mov{X/\P^2}^e$ as
a cone, so by linearity of the action of the Mordell--Weil group,
$K=\R_+ \Pi$ is a fundamental domain for the action of $G$ on
$\Mov{X/\P^2}^e$. In particular $G \cdot K = \Mov{X/\P^2}^e$, which
implies that $\Pic^0(X_\eta) \cdot K = \Mov{X/\P^2}^e$ as required.

\begin{theorem} \label{theorem-funddomain} Suppose $X$ is the blowup
  of $\P^3$ in the base locus of a net of quadrics which has 8
  distinct basepoints, with Mordell--Weil rank $\rho=7$. There
  exists a rational polyhedral cone $U \subset \Mov{X}^e$ such that
  $\Pic^0(X_\eta) \cdot U = \Mov{X}^e$.
\end{theorem} {\bf Proof:} The first step is to choose a rational
polyhedral cone $K_0$ in $\Mov{X}^e$ which maps onto the cone $K$ from
the previous lemma. To do this, we observe that the extremal rays of
$K$ are spanned by the vertices of the polyhedron $\Pi$, which are
vectors of the form $\sum_{i \in I} [E_i] - (|I|-1) [E_1]$, for $I$
any subset of $\{2,\ldots,8\}$. Denote these vectors by
$[v_1],\ldots,[v_n]$. For each vector $[v_i]$ we choose a preimage
$v_i$ in $N^1(X)$ as follows. Suppose $[v_i] = \sum_{i \in I} [E_i] -
(|I|-1) [E_1]$, and put $w_i = \sum_{i \in I} E_i - (|I|-1) E_1$. For
each $[v_i]$, there is an element $\psi_i \in \Pic^0(X_\eta)$ such
that $\psi_i([v_i])=[E_1]$. Therefore in $N^1(X)$, we have
$\psi_i(w_i)=E_1+\frac{m_i}{2}K_X$ for some integer $m_i$. Since the
class $-\frac12K_X$ is preserved by $\Pic^0(X_\eta)$, this gives
$\psi_i(w_i - \frac{m_i+1}{2}K_X)=E_1+\frac{1}{2}K_X$. We then define
$v_i$ to be $w_i - \frac{m+1}{2}K_X$ for $m_i$ chosen as above. The
point of this definition is that each vector $v_i$ then belongs to
$\Mov{X}^e$. For this cone is preserved by pseudo-automorphisms, so it
is enough to show that $\psi_i(v_i) = E_1+\frac{1}{2}K_X$ belongs to
the cone. But this is straightforward: our calculation of $\Curv{X}$
in the previous section shows that $E_1+\frac{1}{2}K_X$ is nef and
hence semi-ample on $X$, so in particular it belongs to $\Mov{X}^e$.
Finally we put $v_0 = -\frac12K_X$ (which also belongs to $\Mov{X}^e$)
and define $K_0$ to be the rational polyhedral cone in $N^1(X)$
spanned by $\{v_0,\ldots,v_n\}$. We have shown that each extremal ray
of $K_0$ lies in the cone $\Mov{X}^e$, so the whole cone $K_0$ does,
and by construction $K_0$ maps onto $K$.

Now choose an SQM $X'$ of $X$ and a divisor $D$ in the ample cone
$A(X')$. By Lemma \ref{lemma_funddomain} there exists a divisor $D_0$
in $K_0$ and an element $\phi \in \Pic^0(X_\eta)$ such that
$\phi_*([D_0]) = [D]$ in $N^1(X/\P^2)$.  Therefore in $N^1(X)$ we have
$\phi_*(D_0)=D+\frac{m}{2}K_X$ for some $m$. If $m \leq 0$ then
$D+\frac{m}{2}K_X$ is also ample on $X'$. If $m>0$ then
$\phi_*(D_0-\frac{m}{2}K_X)=D$ and since $D_0$ belongs to $K_0$ so too
does $D_0-\frac{m}{2}K_X$.  Therefore the union of all translates of
$K_0$ by elements of $\Pic^0(X_\eta)$ intersects the interior of every
nef cone inside $\Mov{X}^e$.

Finally $K_0$ is a rational polyhedral cone in $\Mov{X}^e$ so is
contained in the union $U$ of finitely many nef cones $\Nef{X'}$, each
of which is rational polyhedral by Theorem
\ref{theorem-nefcones}. Since by Lemma \ref{lemma-permutation}
pseudo-automorphisms permute the nef cones of small modifications of
$X$, the union of all the translates of $U$ by elements of
$\Pic^0(X_\eta)$ is a union of nef cones. By the last paragraph this
union intersects the interior of every nef cone, hence equals the
whole effective movable cone $\Mov{X}^e$. \quad QED

\begin{corollary} Suppose $X$ is the blowup
  of $\P^3$ in the base locus of a net of quadrics which has 8
  distinct basepoints, with Mordell--Weil rank $\rho=7$. Then there is
  a finite polyhedral fundamental domain for the action of
  $\PsAut^*(X,\Delta)$ or $\PsAut^*(X)$ on $\Mov{X}^e$.
\end{corollary} {\bf Proof:} If $U$ is the rational polyhedral cone of Theorem
\ref{theorem-funddomain}, we saw that $\Mov{X}^e= \Pic^0(X_\eta) \cdot
U$, hence $\Mov{X}^e=G \cdot U$ for $G$ equal to $\PsAut^*(X,\Delta)$
or $\PsAut^*(X)$. Taking $C$ to be the interior of $\Mov{X}$, this
implies in particular that $C \subset G \cdot U$. Theorem
\ref{theorem-looijenga} then says that $C_+ = G \cdot U = \Mov{X}^e$,
and there is a rational polyhedral fundamental domain for the action
of $G$ on $C_+$.  \quad QED

This completes the proof of the second statement in Theorem
\ref{theorem-mainthm}.

\small
\sc DPMMS, Wilberforce Road, Cambridge CB3 0WB, United Kingdom

Leibniz Universit\"at Hannover, Institut f\"ur Algebraische Geometrie,
Welfengarten 1, D-30167 Hannover, Germany

{\it Email address:} {\tt artie@math.uni-hannover.de}
\end{document}